\documentclass[10pt,a4paper]{article}
\usepackage[T2A]{fontenc}
\usepackage[english,russian]{babel}

\fontsize{12}{21} \selectfont

\newtheorem{theorem}{Theorem}[section]

\newtheorem{definition}{Definition}[section]

\begin{document}

\noindent
\bigskip
\begin{center}
{\large\bf
Another property of minimal surfaces in $E^3$
}
\end{center}
\medskip

\begin{center}
{\bf Andrei~I.~Bodrenko}
\footnote{ Department of Mathematics, Volgograd State University,
Universitє Prospekt 100, Volgograd, 400062, RUSSIA.
Email: bodrenko@mail.ru}
\end{center}

\begin{center}
{\bf Abstract}\\
\end{center}

{\small
The new property of minimal surfaces is obtained in this article.
We derive the equations of $\Delta-$recurrent and
$\Delta-$harmonic surfaces in $E^3$ and
prove that each minimal surface is $\Delta-$recurrent one with
eigenvalue $\varphi=2K$, where $K$ is Gaussian curvature of surface.
We conclude that each minimal surface in $E^3$ satisfies the equation
 $\Delta b = 2Kb$ and
obtain that each $\Delta-$harmonic minimal surface
is a plane $E^2\subset E^3$ or its part.
}

\section*{Introduction}

Let $F^2$ be a smooth two-dimensional surface
in the three-dimensional Euclidean space $E^3$,
$g$ be the induced Riemmanian metric on $F^2$,
$\nabla$ be the Riemmanian connection on $F^2$, determined by $g$,
$b$ be the second fundamental form,
$\overline\nabla$ be the Van der Varden -- Bortolotti covariant derivative,
$\Delta$ be the Laplas operator.

\begin{definition}
\label{definition 1}.
The second fundamental form $b$ is called harmonic if
 $\Delta b \equiv 0$ on $F^2$.
\end{definition}

\begin{theorem}
\label{theorem1}.
If minimal surface $F^2$ in $E^3$ has harmonic
second fundamental form then $F^2$ is a plane
$E^2\subset E^3$ or its part.
\end{theorem}
{\bf Example 1.}
Straight round cylinder $F^2$ is
a $\Delta-$harmonic surface in $E^3$ (see [1]),
i.e. $\Delta b \equiv 0$ on $F^2$.

\begin{definition}
\label{definition 2}.
Surface $F^2$ in $E^3$ is called
$\Delta-$recurrent with eigenvalue $\varphi$,
if the function $\varphi$ on $F^2$ satisfies the condition
$$
\Delta b=\varphi b.
$$
\end{definition}

\begin{definition}
\label{definition 3}.
Surface $F^2$ in $E^3$ is called
$\Delta-$harmonic,
if the second fundamental form $b$ satisfies the condition
$$
\Delta b=0.
$$
\end{definition}

\begin{theorem}
\label{theorem3}.
Each minimal surface $F^2$
 in $E^3$ is $\Delta-$recurrent
with eigenvalue $\varphi=2K$,
where $K$ is the Gaussian curvature of surface.
\end{theorem}

\section{Equations of $\Delta-$recurrent and $\Delta-$harmonic
surfaces in $E^3$}

Let $x$ be an arbitrary point of $F^2$,
$(x^1, x^2, x^3)$ be the Cartesian coordinates
in $E^3$, $(u^1,u^2)$ be the local coordinates
on $F^2$ in some neighborhood $U(x)$ of the point $x$.
Then $F^2$ is given locally by the vector equation
$$
\vec r= \vec r(u^1,u^2) = \{x^1(u^1, u^2), x^2(u^1, u^2), x^3(u^1, u^2)\},
$$
where
$$
\mbox{rg}
\left\|
\frac{\partial x^a}{\partial u^i}
\right\|
= 2\quad , \forall y\in U(x).
$$

Let us fix a point $x\in F^2$ and
 introduce the isothermal coordinates $(u^1, u^2)$
 in some nighborhood $U(x)$ on $F^2$. Then
 the induced metric is $g = A(u^1, u^2)((du^1)^2 + (du^2)^2)$.

Then
$$
g^{11}=g^{22}=\frac{1}{A},\quad g^{12}=0.
$$
We derive the Christoffel symbols:
$$
\Gamma_{11}^1= \Gamma_{12}^2 = -\Gamma_{22}^1 =
\frac{1}{2A}\partial_1 A,
\quad
\Gamma_{12}^1= \Gamma_{22}^2 = -\Gamma_{11}^2 =
\frac{1}{2A}\partial_2 A.
$$
Putting $B=(\ln A)/2$, we have:
$$
\Gamma_{11}^1= \Gamma_{12}^2 = -\Gamma_{22}^1 =
\partial_1 B,
\quad
\Gamma_{12}^1= \Gamma_{22}^2 = -\Gamma_{11}^2 =
\partial_2 B.
$$
We denote by $\{\vec n \}$ the field of unit normal vectors
in the normal bundle $T^\bot F^2$ in $U(x)$:
$\vec n= \vec n(u^1,   u^2)$,
$< \vec n, \vec n> = 1$, where $< , >$ is scalar product in $E^3$.
Let $b_{ij}= <\partial^{2}_{ij}\vec r, \vec n>$ be
the coefficients of the second fundamental form
$b=<d^2\vec r, \vec n> = b_{ij} du^i du^j$.
The covariant derivatives $\nabla_i b_{jk}$ are
$$
\nabla_i b _{jk}=
\partial_i b _{jk} - \Gamma^{m}_{ij} b _{mk} - \Gamma^{m}_{ik} b _{jm}.
$$
The Laplas operator components $\Delta b$ in $U(x)$
are (see [2]):
$$(\Delta b) _{ij}=
g^{kl}\overline{\nabla}_{k}\overline{\nabla}_{l} b_{ij}.
$$
Since the normal connection of $F^2$ in $E^3$
is planar, we have in coordinates $(u^1, u^2)$:
$$A(\Delta b) _{ij}=
\nabla_1\nabla_1 b_{ij}+
\nabla_2\nabla_2 b_{ij}.
$$
We obtain $\nabla_i b _{jk}$:
$$
\nabla_{1}b_{11}=
\partial_1 b_{11}
-2 b_{11} \partial_1 B  + 2 b_{12} \partial_2 B,
$$
$$
\nabla_{1}b_{12}=
\partial_1 b_{12} - 2 b_{12} \partial_1 B + b_{22} \partial_2 B
- b_{11} \partial_2 B,
$$
$$
\nabla_{1}b_{22}=
\partial_1 b_{22} - 2 b_{12} \partial_2 B  - 2 b_{22} \partial_1 B,
$$
$$
\nabla_{2}b_{11}=
\partial_2 b_{11} - 2 b_{11}\partial_2 B  - 2 b_{12} \partial_1 B,
$$
$$
\nabla_{2}b_{12}=
\partial_2 b_{12} - 2 b_{12} \partial_2 B  - b_{22} \partial_1 B
+ b_{11} \partial_1 B,
$$
$$
\nabla_{2}b_{22}=
\partial_2 b_{22} + 2 b_{12} \partial_1 B
- 2 b_{22}\partial_2 B.
$$
The Gauss equation is:
$$
b_{11}b_{22} - b_{12}^2 = -\frac{A}{2}\Delta\ln A.
\eqno (1)
$$
The Peterson-Codacci equation system, written in
 coordinates $(u^1, u^2)$, is:
$$
\partial_2 b_{11} - \partial_1 b_{12} =
(b_{11} + b_{22})\partial_2 B,
\eqno (2)
$$
$$
\partial_1 b_{22} - \partial_2 b_{12} =
(b_{11} + b_{22})\partial_1 B.
\eqno (3)
$$

Let us denote
$\Delta^{e} b_{ij}=\partial^2_{11} b_{ij} +\partial^2_{22} b_{ij}$,
$\Delta B=\partial^2_{11} B +\partial^2_{22} B$.
Using the equations (2) and (3) we have:
$$
A(\Delta b)_{11} = \Delta^{e} b_{11}- 2b_{11}\Delta B -
4(\partial_1 b_{11} +\partial_1 b_{22})\partial_1 B
+
2(b_{11}+b_{22})(3(\partial_1 B)^{2}-(\partial_2 B)^{2}),
$$
$$
A(\Delta b)_{12} = \Delta^{e} b_{12}- 2b_{12}\Delta B
-2(\partial_1 b_{11} +\partial_1 b_{22})\partial_2 B
-
2(\partial_2 b_{11} +\partial_2 b_{22})\partial_1 B
+
$$
$$
+
8(b_{11}+b_{22})\partial_1 B \partial_2 B,
$$
$$
A(\Delta b)_{22} = \Delta^{e} b_{22}- 2b_{22}\Delta B -
4(\partial_2 b_{11} +\partial_2 b_{22})\partial_2 B
+
2(b_{11}+b_{22})(3(\partial_2 B)^{2}-(\partial_1 B)^{2}).
$$
Putting $u=b_{11}+b_{22}$, we obtain:
$$
A(\Delta b)_{11}=\Delta^{e} b_{11}- 2 b_{11} \Delta B -
4 \partial_1 B \partial_1 u + 2u(3(\partial_1 B)^{2}-(\partial_2 B)^{2}),
$$
$$
A(\Delta b)_{12}=\Delta^{e} b_{12}- 2b_{12}\Delta B
- 2\partial_1 B \partial_2 u
- 2\partial_2 B \partial_1 u + 8u\partial_1 B\partial_2 B,
$$
$$
A(\Delta b)_{22}=\Delta^{e} b_{22}- 2b_{22}\Delta B -
4 \partial_2 B \partial_2 u + 2u(3(\partial_2 B)^{2}-(\partial_1 B)^{2}).
\eqno(4)
$$
The surface $F^2$ in $E^3$, by the definition, is
$\Delta-$recurrent with eigenvalue $\varphi$
if the function $\varphi$ satisfies the condition
$$
(\Delta b)_{ij}= \varphi b_{ij}, \quad i, j = 1, 2,
$$
which is equal to the following equation system:
$$
\Delta^{e} b_{11}- 2 b_{11} \Delta B -
4 \partial_1 B \partial_1 u + 2u(3(\partial_1 B)^{2}-(\partial_2 B)^{2})
=A\varphi b_{11},
$$
$$
\Delta^{e} b_{12}- 2b_{12}\Delta B
- 2\partial_1 B \partial_2 u
- 2\partial_2 B \partial_1 u + 8u\partial_1 B\partial_2 B
=A\varphi b_{12},
$$
$$
\Delta^{e} b_{22}- 2b_{22}\Delta B -
4 \partial_2 B \partial_2 u + 2u(3(\partial_2 B)^{2}-(\partial_1 B)^{2})
=A\varphi b_{22}.
\eqno (5)
$$
The equation system (1) --- (3), (5)
determines $\Delta-$recurrent surfaces $F^2$ in $E^3$.

The equality $\Delta b\equiv 0$  on $F^2$ in $E^3$
is equivalent to the following equation system:
$$
1)\quad
\Delta^{e} b_{11}- 2 b_{11} \Delta B -
4 \partial_1 B \partial_1 u + 2u(3(\partial_1 B)^{2}-(\partial_2 B)^{2})
=0,
$$
$$
2)
\quad
\Delta^{e} b_{12}- 2b_{12}\Delta B
- 2\partial_1 B \partial_2 u
- 2\partial_2 B \partial_1 u + 8u\partial_1 B\partial_2 B
=0,
$$
$$
3)
\quad
\Delta^{e} b_{22}- 2b_{22}\Delta B -
4 \partial_2 B \partial_2 u + 2u(3(\partial_2 B)^{2}-(\partial_1 B)^{2})
=0.
\eqno (6)
$$
The equation system (1) --- (3), (6) determines $\Delta-$harmonic surfaces
 $F^2$ in $E^3$.

\section{ $\Delta-$harmonic minimal surfaces}

Proof of the theorem 1. Using the condition on the mean curvature
$H\equiv 0$, we have $u=b_{11}+b_{22}=0$.
Then the equation system (2) and (3) is:
$$
\partial_1 b_{12} = \partial_2 b_{11},
\quad
\partial_2 b_{12} = \partial_1 b_{22}.
$$
Therefore,
$$
\partial^2_{22} b_{11} = \partial^2_{11} b_{22},
\quad
\partial^2_{11} b_{12} + \partial^2_{22} b_{12} =
\partial^2_{21} b_{11} + \partial^2_{12} b_{22}.
\eqno(7)
$$
Using (7), we derive:
$$
\Delta^e b_{11}= \partial^2_{11} b_{11} +\partial^2_{22} b_{11}
=
\partial^2_{11} b_{11} +\partial^2_{11} b_{22} =\partial^2_{11} u = 0,
$$
$$
\Delta^e b_{12}= \partial^2_{11} b_{12} +\partial^2_{22} b_{12}
=
\partial^2_{21} b_{11} +\partial^2_{12} b_{22} =\partial^2_{12} u = 0,
$$
$$
\Delta^e b_{22}= \partial^2_{11} b_{22} +\partial^2_{22} b_{22}
=
\partial^2_{22} b_{11} +\partial^2_{22} b_{22} =\partial^2_{22} u = 0.
\eqno(8)
$$
Using (8), we obtain from the equation system (6):
$$
b_{11} \Delta B = 0,
\quad
b_{12} \Delta B = 0,
\quad
b_{22} \Delta B = 0.
$$
Hence,
$$
(b_{11} b_{22} - b^2_{12})(\Delta B)^2 = 0.
\eqno(9)
$$
Using the equation (1), we obtain from the equation (9):
$$
\left(
-\frac{A}{2}\Delta\ln A
\right)
\left(
\Delta B
\right)^2 = 0.
$$
Consequently, $(\Delta B)^3=0$.
Hence, $\Delta B = 0$ in $U(x)$ and
$$
K=-\frac{1}{2A}\Delta(\ln A)=-\frac{\Delta B}{A} =0.
$$
Therefore, the equations
$$ K = 0, \quad H= 0 $$
are valid in $U(x)$ on $F^2$.

Hence, $U(x)$ is an open part of a plane $E^2\subset E^3$.

The theorem 1 is proved.

\section{$\Delta-$recurrence of
the second fundamental form of\\ minimal surfaces}

Proof of the theorem 2.
From the condition $u=b_{11}+b_{22}=0$ we have, that
the equation system (8) is valid in $U(x)$. Using (8),
 we have from (4):
$$
A(\Delta b)_{11}= - 2 b_{11} \Delta B,
\quad
A(\Delta b)_{12}= - 2b_{12}\Delta B,
\quad
A(\Delta b)_{22}= - 2b_{22}\Delta B.
$$
Therefore, observing
$-\Delta B = A K$,
we derive:
$$
(\Delta b)_{11}=  2K b_{11},
\quad
(\Delta b)_{12}=  2K b_{12},
\quad
(\Delta b)_{22}=  2K b_{22}.
$$
Let us put $\varphi = 2K$.

Therefore,  the equation
$\Delta b =  \varphi b$, where $\varphi = 2K,$
 is valid in $U(x)\subset F^2.$

The theorem 2 is proved.

\bigskip

{\bf Example 2.}
Let the surface $F^2$ is given locally by the vector equation:\\

$\vec
r=\{x^1-\frac{4}{3}(x^1)^3+4(x^1)(x^2)^2,x^2-\frac{4}{3}(x^2)^3+4(x^1)^2(x^2),
2(x^1)^2-2(x^2)^2\}$

Hence,
$g_{11}=(4(x^1)^2+4(x^2)^2+1)^2,$
$g_{22}=g_{11}, g_{12}=g_{21}=0,$ \\
$g^{11}=\frac{1}{g_{11}}, g^{22}=g^{11}, g^{12}=g^{21}=0.$

The unit normal vector:\\
$\vec n= \{ \frac{-16(x^1)(x^2)^2-4(x^1)-16(x^1)^3}{g_{11}} ,
\frac{4(x^2)+16(x^1)^2(x^2)+16(x^2)^3}{g_{11}},
\frac{1-32(x^1)^2(x^2)^2-16(x^1)^4-16(x^2)^4}{g_{11}}\}$

$b_{11}=4, b_{22}=-4, b_{12}=b_{21}=0.$   \\
Therefore, the mean curvature $H=0$
and the Gaussian curvature $ K=\frac{-16}{(g_{11})^2}.$ \\

The covariant derivatives: \\
$\nabla_{1}b_{11}=-\frac{64(x^1)}{\sqrt{g_{11}}},$
$\nabla_{2}b_{12}=\nabla_{2}b_{21}=\nabla_{1}b_{22}=-\nabla_{1}b_{11}$ \\
$\nabla_{2}b_{11}=-\frac{64(x^2)}{\sqrt{g_{11}}},$
$\nabla_{1}b_{21}=\nabla_{1}b_{12}=-\nabla_{2}b_{22}=\nabla_{2}b_{11}$ \\

$\nabla_{2}\nabla_{1} b_{11}=\frac{3584(x^1)(x^2)}{g_{11}},$ \\
$\nabla_{1}\nabla_{2} b_{11}=\nabla_{1}\nabla_{1} b_{12}=
\nabla_{2}\nabla_{2} b_{12}=\nabla_{1}\nabla_{1} b_{21}=\\
-\nabla_{2}\nabla_{2} b_{21}=-\nabla_{2}\nabla_{1} b_{22}=
-\nabla_{1}\nabla_{2} b_{22}=\nabla_{2}\nabla_{1} b_{11}.$

$\nabla_{1}\nabla_{1} b_{11}=\frac{64(-28(x^2)^2-1+28(x^1)^2)}{g_{11}},$ \\
$\nabla_{1}\nabla_{2} b_{12}=\nabla_{1}\nabla_{2} b_{21}=
\nabla_{1}\nabla_{1} b_{22}=-\nabla_{1}\nabla_{1} b_{11}.$ \\

$\nabla_{2}\nabla_{2} b_{11}=-\frac{64(-28(x^2)^2+1+28(x^1)^2)}{g_{11}},$ \\
$\nabla_{2}\nabla_{1} b_{12}=\nabla_{2}\nabla_{1} b_{21}=
-\nabla_{2}\nabla_{2} b_{22}=\nabla_{2}\nabla_{2} b_{11}.$

$(\Delta b)_{11} = -\frac{128}{(g_{11})^2},
(\Delta b)_{22}=-(\Delta b)_{11},
(\Delta b)_{12}=0,
(\Delta b)_{21}=0.$

Therefore, we have $\Delta b = 2Kb.$

{\bf Example 3.}
Let the surface $F^2$ is given locally by the vector equation:\\

$\vec r=\{x^2\cos(x^1), x^2\sin(x^1), x^1\} $    \\

$g_{11}=1+(x^2)^2, g_{22}=1, g_{12}=g_{21}=0. $ \\

$g^{11}=\frac{1}{1+(x^2)^2}, g^{22}=1,g^{12}=g^{21}=0. $ \\

The unit normal vector:\\
$\vec n=\{-\frac{\sin(x^1)}{\sqrt{1+(x^2)^2}}, \frac{\cos(x^1)}{\sqrt{1+(x^2)^2}}, \frac{-x^2}{\sqrt{1+(x^2)^2}}\} $ \\

$b_{11}=b_{22}=0, b_{12}=b_{21}=\frac{1}{\sqrt{1+(x^2)^2}}. $  \\

Therefore, the mean curvature $H=0$
and the Gaussian curvature $K=-\frac{1}{(1+(x^2)^2)^2}.$

The covariant derivatives:\\
$\nabla_{1}b_{11} = \frac{2(x^2)}{\sqrt{1+(x^2)^2}},\\
\nabla_{2}b_{11}=\nabla_{1}b_{12}=\nabla_{1}b_{21}=\nabla_{2}b_{22}=0 \\
\nabla_{2}b_{12}=\nabla_{2}b_{21}=\nabla_{1}b_{22}=
\frac{-2x^2}{\sqrt{(1+(x^2)^2)^3}}. $

$\nabla_{1}\nabla_{1} b_{11}=\nabla_{2}\nabla_{2} b_{11}=
\nabla_{2}\nabla_{1} b_{12}=\nabla_{1}\nabla_{2} b_{12}=
\nabla_{2}\nabla_{1} b_{21}=\nabla_{1}\nabla_{2} b_{21}=
\nabla_{1}\nabla_{1} b_{22}=\nabla_{2}\nabla_{2} b_{22}=0.$ \\

$\nabla_{1}\nabla_{2} b_{11}=\nabla_{1}\nabla_{1} b_{12}=
\nabla_{1}\nabla_{1} b_{21}=-\nabla_{1}\nabla_{2} b_{22}=
-\frac{6(x^2)^2}{\sqrt{(1+(x^2)^2)^3}}. $ \\

$\nabla_{2}\nabla_{1} b_{11}=
-\frac{2(-1+3(x^2)^2)}{\sqrt{(1+(x^2)^2)^3}}$ \\

$\nabla_{2}\nabla_{2} b_{12}=
\nabla_{2}\nabla_{2} b_{21}=\nabla_{2}\nabla_{1} b_{22}=
\frac{2(-1+3(x^2)^2)}{\sqrt{(1+(x^2)^2)^5}}.$ \\

$(\Delta b)_{12} = -\frac{2}{\sqrt{(1+(x^2)^2)^5}},
(\Delta b)_{21}=(\Delta b)_{12},
(\Delta b)_{11}=0,
(\Delta b)_{22}=0.$ \\

Therefore, we have $\Delta b = 2Kb.$ \\

{\bf References} \\
\begin{enumerate}
\item Bodrenko A.I.,
Developable laplas recurrent curves in Euclidean space $E^{3}$.
Surveys in Applied and Industrial Mathematics,
2001, Vol. 8, Issue 2,  p. 540.
\item Bodrenko A.I.,
Laplas recurrent curves in Euclidean spaces.
//Vestnik of Volgograd State University.
Vol. 9. Issue 1. 2001. Part 2. pp. 6--8.
\end{enumerate}

\end{document}